
\def\A{{\cal A}}
\def\B{{\cal B}}
\def\C{{\cal C}}
\def\P{{\cal P}}
\def\L{{\cal L}}
\def\M{{\cal M}}
\def\V{{\cal V}}
\baselineskip=14pt
\parskip=10pt
\def\halmos{\hbox{\vrule height0.15cm width0.01cm\vbox{\hrule height
  0.01cm width0.2cm \vskip0.15cm \hrule height 0.01cm width0.2cm}\vrule
  height0.15cm width 0.01cm}}
\font\eightrm=cmr8 
\font\eighttt=cmtt8
\magnification=\magstephalf
\def\P{{\cal P}}

\def\1{{\overline{1}}}
\def\2{{\overline{2}}}
\parindent=0pt
\overfullrule=0in

\def\frac#1#2{{#1 \over #2}}
\bf
\centerline
{
${{5} \choose {2}}$ Proofs that ${{n} \choose {k}} \leq {{n} \choose {k+1}}$ if $k<n/2$.
}
\rm
\bigskip
\centerline{ {\it
Doron 
ZEILBERGER}\footnote{$^1$}
{\eightrm  \raggedright
Department of Mathematics, Rutgers University (New Brunswick),
Hill Center-Busch Campus, 110 Frelinghuysen Rd., Piscataway,
NJ 08854-8019, USA.
{\eighttt zeilberg  at math dot rutgers dot edu} ,
\hfill \break
{\eighttt http://www.math.rutgers.edu/\~{}zeilberg/} .
March 4, 2010.
An almost verbatim transcript of a mathematics Colloquium talk delivered at
Columbia University (Mathematics 520), Feb. 17, 2010, 5:00-6:00pm EST.
I would like to thank Mikhail Khovanov for inviting me, and the
nineteen brave souls
who attended (for example Walter Neumann). 
I also would like to express my great disappointment at the 
71 (faculty)+58 (grad students)-19=110
people (for example David Bayer) who did {\it not} attend.
Supported in part by the NSF.
Exclusively published in the Personal Journal of Ekhad and Zeilberger 
({\eighttt http://www.math.rutgers.edu/\~{}zeilberg/pj.html}), and {\eighttt arxiv.org}.
}
}

There is no trivial mathematics, there are only trivial mathematicians!
A mathematician is {\it trivial} if he or she believes that there exists trivial mathematics.
But this is not the only way for a mathematician to be trivial. Another sufficient
condition for a mathematician to be trivial is {\it not} to show up to a colloquium talk
with such an intriguing title and abstract! Conversely, if you {\it do} show up,
you are definitely non-trivial, so, congratulations, dear audience, you are the
(only!) $19$ non-trivial mathematicians in Columbia University.

There are (at least) {\it two ways} to define ${{n} \choose {k}}$. One way is
$$
{{n} \choose {k}}:={{n!} \over {k!(n-k)!}} \quad .
$$
By this definition, the statement of the title is indeed trivial:
$$
{
{{{n} \choose {k}}}
\over
{{{n} \choose {k+1}}}
}
=
{
{{{n!} \over {k!(n-k)!}}}
\over
{{{n!} \over {(k+1)!(n-k-1)!}}}
}
=
{
{(k+1)! (n-k-1)!}
\over
{k!(n-k)!}
}
=
{{k+1} \over {n-k}} \leq 1 \quad, 
$$
if $k<n/2$.

The other definition is a {\it combinatorial} one.  ${{n} \choose {k}}$ is the number of ways of
{\it choosing} a set of $k$ members out of an $n$-element set. It is also the
number of $n$-letter words in the alphabet $\{Street, Avenue\}$ with exactly $k$ occurrences
of the ``letter'' ``Street'', as I was reminded when I walked, earlier today, in the (real!, not proverbial)
Manhattan lattice from Pennsylvania Station to the Columbia campus (except I confess that I cheated,
and walked most of the way on Broadway). If you adapt the latter definition, then 
${{n} \choose {k}}={{n!} \over {k!(n-k)!}}$ becomes a {\it theorem}, that can be proved, e.g., by proving
that both sides satisfy the recurrence (and initial condition).
$$
f(n,k)=f(n-1,k-1)+f(n-1,k) \quad , \quad
f(0,k)=\delta_{0,k}  \quad .
$$
(If after you walked $n$ blocks, you are currently at the corner of $k$-th Street and $(n-k)$-th Avenue, then one
block earlier, you were either at  the corner of $(k-1)$-th Street and $(n-k)$-th Avenue or
the corner of $k$-th Street and $(n-k-1)$-th Avenue).

In particular, it follows that ${{n!} \over {k!(n-k)!}}$ is always an integer!, which is not so obvious,
(since this is a {\it ratio} of two integers, that 
morally should be a fraction, unless some miracle occurs),
and that it is less than $2^n$. The special case that 
$$
{{(2n)!} \over {n!^2}}={{(2n)(2n-1)\cdots (n+1)} \over {(1) \cdots (n)} } \quad
$$
is an integer, and that it is less than $2^{2n}$, has an {\it enormous} number-theoretical significance.
It was used by Chebychev, in 1851, to ``almost'' prove the Prime Number Theorem. 
Even though this breakthrough was ``superseded'' by the full Prime Number Theorem, first proved
at the end of the 19th-century, all the proofs of 
the latter, as well as the later elementary proofs
of Erd{\H o}s and Selberg, use Chebychev's result as a stepping-stone for the stronger statement.
More recently, it turned out to be {\it crucial} in the 
amazing Agrawal-Kayal-Saxena[AKS] $PRIMES \in \P$ proof.

Let's recall Chebychev's argument. Since people today are 
{\it so} specialized, I am willing to bet
that many of you have never seen it before. Only this gem is worth the 
{\it admission fee} of this talk (which is an hour of your precious time,
that at least \% 80 of the Columbia faculty and graduate students found too
exorbitant.)

Let's look at all the prime numbers between $n$ and $2n$. They must all divide the
{\it integer} ${{(2n)!} \over {n!^2}}$, so
$$
\prod_{n \leq p \leq 2n} p \leq 2^{2n}  \quad,
$$
Now take $log$ of both sides, define $\theta(x)=\sum_{p \leq n} \log p$, and you would get that
$$
\theta(2n)-\theta(n) \leq (2n) \log 2 \quad ,
$$
that implies
$$
\theta(n)-\theta(n/2) \leq (n) \log 2
$$
$$
\theta(n/2)-\theta(n/4) \leq (n/2) \log 2
$$
$$
\dots
$$
Adding these up, you get that $\theta(n) \leq (2 \ln 2) n$, which is equivalent to
$\pi(x) \leq C {{x} \over {\log x}}$ for $C=2 \ln 2=1.386 \dots	$. 
Later Chebychev made $C$ even smaller, and Sylvester got very close to $1$, and
analogously for lower bounds, but the full Prime Number Theorem had to wait for
Hadamard and de la Vall\'ee Poussin, in 1896.

Going back to proving (and reproving) that
${{n} \choose {k}} \leq {{n} \choose {k+1}}$ if $k<n/2$, here is an {\it inductive} proof.
$$
{{n} \choose {k}}= {{n-1} \choose {k-1}} + {{n-1} \choose {k}} 
\leq {{n-1} \choose {k}} + {{n-1} \choose {k+1}}  \, = \, {{n} \choose {k+1}} \quad,
$$
by the {\it induction hypothesis}, (provided that the {\it hypothesis} is fulfilled!). 
This is always true for ${{n-1} \choose {k-1}} \leq {{n-1} \choose {k}}$, since
$k<n/2$ implies $k-1 < (n-1)/2$, but for 
${{n-1} \choose {k}} \leq {{n-1} \choose {k+1}}$ it may happen that
$k<n/2$ but $k \geq (n-1)/2$. This happens exactly when $n=2k+1$, and for this special case
we have to separately prove:
$$
{{2k+1} \choose {k}} \leq {{2k+1} \choose {k+1}} \quad ,
$$
but this follows from the {\it even stronger} fact that:
$$
{{2k+1} \choose {k}} = {{2k+1} \choose {k+1}} \quad ,
$$
by the symmetry of the binomial coefficients.

I admit that this is an {\it ugly duckling} of a proof ({\it manipulatorics}, induction), but
by carefully tracing it, we can get a {\it beautiful swan} of a proof, by
defining an {\it explicit} injection that maps, in a {\it canonical way},
an $n$-letter word in the alphabet $\{S,A\}$ with $k$ S's to one with $k+1$ S's.
Simply look at the last time the number of Avenues exceeded the number of 
Streets by exactly one, and swap Avenues  and Streets,
until then, and leave the rest intact.

The proof that I just gave is an example of a {\bf combinatorial proof}, and the process
of finding a {\it combinatorial interpretation} to an algebraic identity or inequality,
using a bijection and injection respectively, is called {\bf Combinatorization}.
Often algebraic/inductive proofs can be ``traced'' and converted to beautiful bijective or
injective proofs, like in the above case.

Speaking of {\it combinatorization}, this is the grand-daddy of a more recent trend, called
{\it categorification}, made popular by master-blogger John Baez. Categorification became
a {\it household name} when my host, Mikhail Khovanov[Kh], in 2000, astounded the
mathematical world by {\it categorifying} the famous Jones polynomials, by replacing a boring
polynomial by an exciting {\it cell-complex}. I strongly recommend Dror Bar-Natan's ([B]) very lucid
exposition of Khovanov's seminal ideas, that you can easily find in {\tt arxiv.org}.

Going back to combinatorics, we will meet other, even better, combinatorial 
proofs, later on, but let me now present to you
yet another {\it algebraic proof}. 
Using the {\it Zeilberger algorithm}[Z2] (or otherwise\footnote{$^2$}
{\eightrm  \raggedright
Algebra is really combinatorics in disguise, when you expand $(1+x)^n$ you make $n$ independent decisions,
whether to pick the $1$ or the $x$. 
The coefficient of $x^k$ is the the number of ways of
{\it choosing} which $k$ of the $n$ terms will donate its $x$
to the common cause.}),
we can find the
{\it generating function}
$$
P_n(x):=\sum_{k=0}^{n} {{n} \choose {k}} x^k = (1+x)^n \quad,
$$
and another way of stating that
${{n} \choose {k}} \leq {{n} \choose {k+1}}$ if $k<n/2$
is to say that the coefficients of $P_n(x)=(1+x)^n$
first go up and then go down. Such a polynomial is called {\it unimodal}. In our case it is
also {\it symmetric}. Let's call a {\it symmetric and unimodal} polynomial with non-negative
integer coefficient a $Z$-polynomial, and let's call the {\it darga} of $P$ its
low-degree plus its (high)-degree. For example, the {\it darga} of $x^4+x^5$ is $9$ while
the {\it darga} of $x^3$ is $6$.

The following two simple facts (taken from my {\it de-combinatorization} ([Z1])
of Kathy O'Hara's ([O]) seminal combinatorial proof of the unimodality of the Gaussian polynomials) are easily proved.

{\bf Fact 1}: The sum of two Z-polynomials of the same darga is another Z-polynomial of that darga.

{\bf Fact 2}: The product of two Z-polynomials is yet-another-one, and its darga is the sum of their dargas.

To prove Fact 2 note that the additive
``atoms'' of Z-polynomials are polynomials of the form
$$
x^i+x^{i+1}+...+x^{j}
$$
and multiplying out two such atoms would yield
$$
(x^a+x^{a+1}+...+x^{b})(x^{c}+x^{c+1}+...+x^{d}) = x^{a+c}+2x^{a+c+1}+3x^{a+c+2}+ \dots + 3x^{b+d-2} + 2x^{b+d-1} + x^{b+d} \quad ,
$$
which is indeed a Z-polynomial of darga $(a+c)+(b+d)=(a+b)+(c+d)$.

It follows immediately, by induction, that $(1+x)^n$ is a Z-polynomial, since $1+x$ is. But we get, for
the same price, that many other polynomials are Z-polynomials, 
and hence automatically unimodal. For example
$$
(x+x^2+x^3+x^4+x^5+x^6)^n \quad,
$$
which has the following probabilistic interpretation. You roll a fair die $n$ times and at each roll
you  win as many dollars as the number of dots that show up. Then you are more likely to win
$k+1$ dollars than $k$ dollars as long as $k$ is less than your expected gain $7n/2$.

More generally:
$$
(1+x)^m(x+x^2+x^3+x^4+x^5+x^6)^n (x+x^2+x^3+x^4)^k \quad ,
$$
that also has a gambling interpretation, and 
many more complicated gambling scenarios, that you are welcome to make up.

Let's take a closer look at the above combinatorial proof that
${{n} \choose {k}} \leq {{n} \choose {k+1}}$, that consisted 
in defining an explicit injection between $k$-sets to $(k+1)$-sets.
It inputs a set $S$ with $k$ elements ($k<n/2$) and outputs 
a set with one more element by looking at the smallest integer $r$
such that $|S \cap \{1,2, \dots , r \}|=(r-1)/2$ and mapping it
to the set $(\{1, \dots ,r \} \backslash S) \cup (S \cap \{r+1, \dots , n \})$.

For example, with $n=11$ and $k=4$ the $4$-set $\{1,2,4, 11\}$ is mapped to the $5$-set
$\{3,5,6,7, 11\}$ (in this example $r=7$). Note that for this injection
the output-set does {\it not} contain the input set. 
It would be more desirable, and {\it natural}, if we could come-up with an injection
$S \rightarrow S'$ from the collection of $k$-sets to the collection of $k+1$-sets that
has the property that $ S \subset S'$, in other words, find a ``rule'' that
adds a new member to $S$, as long as $k<n/2$, and in such a way that
no two different $S$s would give the same $S'$.

If there would be such a mapping we would get, by iterating it, a {\it maximal chain}
that ends at the middle rank. By symmetry, if we reflect this to the complement,
we would get a {\it central chain decomposition} of the Boolean lattice (alias $n$-dimensional
unit cube). Conversely, any such chain-decomposition of the Boolean lattice
would give such an injection, and would yield yet-another-proof of the unimodality
of the binomial coefficients.

The easiest way to construct such a chain decomposition is {\it recursively}.
Take any symmetric chain of $B_{n-1}$
$$
C_r \rightarrow C_{r+1} \rightarrow \dots \rightarrow C_{n-1-r} \quad,
$$
and construct two new chains in $B_n$. The first is the same
$$
C_r \rightarrow C_{r+1} \rightarrow \dots \rightarrow C_{n-1-r} \quad,
$$
but viewed as belonging to $B_n$, and the second is
$$
C_r \cup {n} \rightarrow C_{r+1} \cup {n} \rightarrow \dots \rightarrow C_{n-1-r} \cup {n}\quad .
$$
There is only one problem! Neither chains are legitimate symmetric chains in $B_n$.
The sum of the starting rank and ending rank (in the case of the Boolean lattice, the
rank of a set is its number of elements) should be $n$, whereas the first chain has
the sum {\it too low}, namely $n-1$, while the second chain has its rank {\it too high}, namely
$n+1$. To get two new chains that are {\it just right}, we cut the last member of the
second chain and put it at the end of the first chain, getting the two chains:
$$
C_r \rightarrow C_{r+1} \rightarrow \dots \rightarrow C_{n-1-r} \rightarrow C_{n-1-r} \cup {n} \quad,
$$
and
$$
C_r \cup {n} \rightarrow C_{r+1} \cup {n} \rightarrow \dots \rightarrow C_{n-2-r} \cup {n}\quad .
$$

Let's illustrate this construction for $n \leq 3$.
For $n=1$ we only have one chain, namely:
$$
\emptyset \rightarrow \{ 1 \} \quad .
$$
This gives rise to two chains for $n=2$:
$$
\emptyset \rightarrow \{ 1 \} \rightarrow \{1,2\} \quad, \quad \{ 2 \} \quad .
$$
The first of these gives rise to two chains for $n=3$:
$$
\emptyset \rightarrow \{ 1 \} \rightarrow \{1,2\} \rightarrow \{1,2,3\} \quad , \quad
\{ 3 \} \rightarrow \{ 1,3 \} \quad ,
$$
while the singleton chain $\{ 2 \}$ only gives rise to {\it one} chain (the second one is empty)
$$
\{ 2 \} \rightarrow \{2,3\} \quad .
$$

Martin Aigner came up with another way of constructing a symmetric chain decomposition
for the Boolean lattice $B_n$, that may be termed {\it lexicographic greed}.
Start with the empty set, and at each level look at the lexicographically first set that has not yet
been committed and that contains the current tail of the emerging chain. Keep doing it until you get stuck.
Surprisingly, you get a symmetric chain decomposition. Why?, because it happens to be the same as
the one above. So even though many people would find Aigner's construction more elegant and
appealing, the easiest way to prove its validity is to discover the recursive construction above
and then it is easy to prove by induction that it is indeed the same.

The drawback that both the recursive and Aigner's([A]) lexicographic-greed approaches share is that you  have
to construct {\it all chains}, and find out how the injection acts on {\it all} sets, at once,
requiring exponential time and space. What if you only care about
the successor of just {\it one} individual set? Curtis Greene and Daniel Kleitman[GK] came up
with a very elegant description of (essentially the same!) injection.

There is a one-to-one mapping between sets and words in the {\it alphabet}
$\{ [ , ] \}$. For any set $S$ of $n$ elements form the ``word''
$(w_1, \dots, w_n)$ by the rule $w_i= [$ iff $i \in S$.
For example, the empty set for $n=4$ corresponds to the word $]]]]$ and the
whole set $\{1,2,3,4\}$ corresponds to the word $[[[[$.

If you have a {\it legal bracketing} then it forms its own singleton chain. Otherwise, ``compile'' it
to the best of your ability, matching a left-bracket ``['' with a right one ``]''. 
Once you have finished ``compiling'' you would get a bunch of ]'s followed by
a bunch of ['s which is as illegal as it gets, possibly (and usually) interspersed
with clusters of legal bracketings. Leave these legal bracketings alone, and
change the last ] by a [. In symbols:
$$
L_1 \quad ] \quad L_2 \quad ] \quad L_3 \quad ] \quad \dots \quad L_{k-1} \quad ]
\quad L_k \quad ] \quad L_{k+1} \quad [ \quad L_{k+2} \quad [ \quad \dots \quad [  \quad L_r \quad
$$
goes to
$$
L_1 \quad ] \quad L_2 \quad ] \quad L_3 \quad ] \quad \dots \quad L_{k-1} \quad ]
\quad L_k \quad [ \quad L_{k+1} \quad [ \quad L_{k+2} \quad [ \quad \dots \quad [  \quad L_r \quad
$$
If you can't do it (i.e. $k=0$), then the chain ends.

The existence of a symmetric chain decomposition for the Boolean lattice immediately implies
{\bf Sperner's theorem} that the largest possible collection of subsets of $\{1,2, \dots n\}$
such that none of its members
properly contains  another one (what is called an {\it anti-chain}, or {\it clutter})
equals ${{n} \choose {[n/2]}}$. Obviously, this is sharp, since the 
collection of {\it all} $[n/2]$-sets , that has ${{n} \choose {[n/2]}}$ members, is obviously {\it an} anti-chain. 
Can you do better? Of course not! 
The number of symmetric chains in any symmetric chain decomposition of the Boolean lattice
(and we know that one exists) equals ${{n} \choose {[n/2]}}$, since each chain passes once through the middle-rank, and
every set belongs to exactly one chain. Given {\it any} anti-chain,
there can be {\it at most} one-set-per-chain, or else it would not be an anti-chain!

While the above proof of Sperner's theorem is my personal favorite,
let me remind you of another, just-as-nice  {\it proof from the book}, due to David Lubell[L].

There are $n!$ possible chains that start at the top, the empty set, and end-up at the
bottom ($\{1,2, \dots n\}$) (in obvious one-one correspondence with permutations).
Let $\C$ be a potential anti-chain. For each $S \in \C$, there are exactly
$|S|!(n-|S|)!$ such top-to-bottom chains that pass through $S$, and of course,
no two different members of $\C$ can share such a top-to-bottom chain, or else they
would be related!

So we have the obvious inequality
$$
\sum_{S \in \C} |S|!(n-|S|)! \, \leq n! \quad,
$$
that implies that
$$
\sum_{S \in \C} {{1} \over 
{ 
{{n} \choose {|S|}}
}
}  
\, \leq 1 \quad .
$$
But the maximum of ${{n} \choose {|S|}}$ is ${{n} \choose {[n/2]}}$
(thanks to the main theorem of the present article!), so the minimum of
$1/{{n} \choose {|S|}}$ is $1/{{n} \choose {[n/2]}}$, and we have
$$
{{|\C|} \over {{{n} \choose {[n/2]}}}} \, \leq \, 1 \quad ,
$$
as claimed.

{\bf The Last (and longest! (yet the best!))  Proof}

We have already presented above several {\it combinatorial} proofs of
$$
{{n} \choose {k}} \leq {{n} \choose {k+1}} \quad, \quad if \quad k<n/2 \quad ,
$$
by finding a {\it set-theoretical} {\bf injection} between the collection of $k$-sets and
the collection of $(k+1)$-sets, i.e. between two {\it sets} (of sets).

In general, a {\it combinatorial} proof of 
$$
a \leq b
$$
consists of constructing {\it sets} $A$ and $B$ such $a=|A|$ and $b=|B|$, and an injection
$$
f: A \rightarrow B \quad .
$$

But, there is {\it yet another way}, a linear-algebra proof! Come-up with two {\it vector spaces}
$\A$ and $\B$ such that $dim(\A)=a$ and $dim(\B)=b$
and construct a {\it linear transformation}
$$
T: \A \rightarrow \B \quad ,
$$
and prove that $T$ is an injection by proving that for any $f \in \A$, $Tf=0$ implies $f=0$.

Let $\V_k$ be the vector space spanned by all $k$-subsets of $\{1, \dots , n \}$, in other words
the vector space of all ``formal sums'' (as they would say in algebraic topology)
$$
\sum_{|S|=k} a_S S \quad ,
$$
where $a_S$ are members of your favorite field (say the field of rational numbers, or even $GF(p)$ for any prime $p$ larger than $n$).

Our proposed mapping, $\M: \V_k \rightarrow \V_{k+1}$,
soon to be proved an injection, is defined on basis elements by
$$
\M (S)= \sum_{j \not \in S} (S \cup {j}) \quad ,
$$

and extended linearly. What is the ``meaning'' of $\M(S)$? 
Suppose that you enlarge your current faculty $S$ by another member, and you can't decide,
and you want to hire everyone who is not already in $S$ but you are only allowed to hire one person.
If you live in a classical world, you would have to make-up your mind, make one new professor happy, but
disappoint all the other applicants. But in the quantum world, you can have a ``superposition'' of all
scenarios for ``hiring an extra professor''.

In order to prove that $\M$ is indeed an injection, we need a ``companion operator'': 
$\L: \V_k \rightarrow \V_{k-1}$, defined on basis elements by:
$$
\L (S)= \sum_{i \in S} (S \backslash {i}) \quad ,
$$
and extended linearly. $\L(S)$ has an analogous meaning in a quantum world. Because of budget cuts, you have
to fire one professor, but you don't want to get anyone upset, so you have a quantum-superposition of
all firing-one-professor scenarios.

I now claim that on $\V_k$,
$$
\M \L - \L \M = \mu(k) I \quad,
$$
where $\mu(k)=2k-n$ ($n$ is fixed through this proof), and $I$ is the identity mapping. 
Of course, by linearity, it is enough to prove this for basis elements $S \in \V_k$:
$$
\M \L (S) - \L \M (S)  = \mu(k) S \quad . 
\eqno(1)
$$
$\M\L(S)$ is formal sum of all scenarios of fire-and-then-hire while $\L\M(S)$ is the formal sum
of all scenarios of hire-and-then-fire. If the guy you hired and the guy you fired are {\it different} 
then ``hire-Smith-then-fire-Jones'' yields the same set as ``fire-Jones-then-hire-Smith'' and so they
cancel out. The only scenarios that do not cancel out are 
those where the guy you fired and the guy you hired are
one and the same. There are $k$ ways to fire-and-then-hire the same person, and there are $n-k$ ways
to hire-and-then-fire the same person, at each case resulting in the original set $S$.
This gives a net contribution of $k-(n-k)=2k-n$ copies of $S$.

Next I claim that on $\V_k$, for any $r \geq 1$
$$
\M \L^r - \L^r \M = (\mu(k)+ \dots +\mu(k-r+1)) \L^{r-1} \quad .
\eqno(2)
$$
This follows easily by induction on $r$, by using
$$
\M \L^{r+1} - \L^{r+1} \M =  (\M \L^{r} - \L^{r} \M ) \L+ \L^r( \M \L - \L \M) \quad .
$$
So, if $f \in \V_k$, we have:
$$
( \M \L^{r+1} - \L^{r+1} \M )f =  ( \M \L^r - \L^r \M ) ( \L f) + \L^r ( \M \L - \L \M) f \quad .
\eqno(3)
$$
Since $\L f \in \V_{k-1}$, we have from the induction hypothesis that the first term on the right side of $(3)$ equals
$$
(\M \L^{r} - \L^{r} \M ) (\L f)= (\mu(k-1)+ \dots + \mu(k-r)) \L^{r-1} (\L f) = (\mu(k-1)+ \dots + \mu(k-r)) \L^r f \quad ,
\eqno(3a)
$$
and since $( \M \L - \L \M) f=\mu(k)f$, the second term of $(3)$ is
$$
\L^r ( \M \L - \L \M) f= \L^r \mu(k) f= \mu(k) \L^r f \quad .
\eqno(3b)
$$
Incorporating $(3a)$ and $(3b)$ into $(3)$, we get:
$$
( \M \L^{r+1} - \L^{r+1} \M )f =  ( \M \L^r - \L^r \M ) ( \L f) + \L^r ( \M \L - \L \M) f=
$$
$$
(\mu(k-1)+ \dots + \mu(k-r)) \L^r f + \mu(k) \L^r f =
(\mu(k)+ \dots + \mu(k-r)) \L^r f  \quad ,
$$
that is $(2)$ with $r$ replaced by $r+1$.

Now suppose that there is an $f \in \V_k$ such that $\M f=0$. We have to prove that $f=0$.
By $(2)$ we have that
$$
\M \L^r f =(\mu(k)+ \dots +\mu(k-r+1)) \L^{r-1} f \quad .
$$
Applying $\M^{r-1}$ to both sides gives
$$
\M^r \L^r f =(\mu(k)+ \dots +\mu(k-r+1)) \M^{r-1}\L^{r-1} f \quad .
$$
Iterating, gives:
$$
\M^r \L^r f =(\mu(k)+ \dots +\mu(k-r+1)) (\mu(k)+ \dots +\mu(k-r+2))  \cdots (\mu(k)) f \quad .
$$
So we have
$$
\M^{k+1} \L^{k+1} f= (Non-Zero-Number) f \quad .
$$
But, since $f \in \V_k$, $\L^{k} f$ is a multiple of the empty set, and hence $\L^{k+1} f=0$
($\L \emptyset =0$, since in that case we get the {\it empty sum} in the definition of $\L \emptyset$).
So we get that $f=0$, as promised.

So indeed, if $k < n/2$, the mapping $\M: \V_k \rightarrow \V_{k+1}$ is an injection, and we get
$dim(\V_k) \leq dim(\V_{k+1})$, and so, once again, we know that
${{n} \choose {k}} \leq {{n} \choose {k+1}}$ if $n<k/2$. \halmos

But why work so hard, if we had the former far easier proofs? 
One reason, is {\it why not}? Who said that an elegant proof has to be short?
Another reason is that this proof {\it extends}, almost verbatim, to other lattices,
for which no simple proofs of rank-unimodality and the Sperner property are known.
The proof that I just presented was inspired by, and 
is along similar lines as-but not quite the same-
as Robert Proctor's ([P]) beautiful simplification of Richard Stanley's([S])  seminal proof
of the Sperner property for lattices of integer partitions.
The main part in the proof of Spernerity,
proving that $\M$ is injective (as we just did), can be traced, in an almost
equivalent form 
(but using differential operators operating on so-called semi-invariants)
to James Joseph Sylvester[Sy] way back in 1878.

{\bf References}

[A] Martin Aigner, {\it Lexicographic matching in Boolean algebras}, J. Combinatorial Theory Ser. B {\bf 14} (1973), 187--194. 

[AKS] Manindra Agrawal, Neeraj Kayal and Nitin Saxena, {\it PRIMES is in P}, Annals of Mathematics, 
Second Series {\bf 160} (2004), 781-793 .

[B] Dror Bar-Natan, {\it On Khovanov's categorification of the Jones polynomial}, Algebr. Geom. Topol. (electronic) {\bf 2}  (2002), 337--370 .

[GK] Curtis Greene and Daniel J. Kleitman, {\it Proof techniques in the theory of finite sets},
Studies in combinatorics,  MAA Stud. Math. {\bf 17}, (G. -C. Rota, ed.), Math. Assoc. America, Washington, D.C., 1978

[K] Mikhail Khovanov, {\it A categorification of the Jones polynomial},  Duke Math. J.  {\bf 101}  (2000),  359--426 .

[L] David Lubell, {\it A short proof of Sperner's lemma},  J. Combinatorial Theory  {\bf 1}  (1966),  299.

[O] Kathleen M. O'Hara, {\it Unimodality of Gaussian coefficients: a constructive proof},
 J. Combin. Theory Ser. A  {\bf 53} (1990),  29--52. 

[P] Robert A. Proctor,
{\it Solution of two difficult combinatorial problems with linear algebra}, Amer. Math. Monthly  {\bf 89}(1982), 721--734

[St] Richard P. Stanley, {\it Weyl groups, the hard Lefschetz theorem, and the Sperner property},
SIAM J. Algebraic Discrete Methods  {\bf 1}(1980), 168--184.

[Sy] James Joseph Sylvester, {\it Proof of the hitherto undemonstrated fundamental theorem of invariants},
Philosophical Magazine {\bf 5} (1878), 178-188. (Also in Collected Works v. 3, 117-126, Chelsea, NY, 1973).

[Z1] Doron Zeilberger, {\it A one-line high school algebra proof of the unimodality of the Gaussian polynomials $[^n_k]$ for $k<20$}.  
In: $q$-series and partitions (Minneapolis, MN, 1988),  D. Stanton, ed., 67--72, IMA Vol. Math. Appl., {\bf 18} (1989), Springer, New York.

[Z2] Doron Zeilberger, {\it The method of creative telescoping}, J. Symbolic Comput. {\bf 11} (1991), 195--204. 

\end